\newcommand{\abstand}{\vspace{1em}}

\newtheorem{theo}{Theorem}
\newtheorem{remark}{Remark}

\newcommand{\proof}{{\em Proof. }}
\newcommand{\qed}{$\hfill\Box$}

\newcommand{\zitat}[4]{\bibitem{#1}{\sc #2}, {\em #3\/}. #4.\vspace{-0.0em}}

\newcommand{\PG}[2]{\mbox{$\mbox{{\rm PG}}(#1,#2)$}}
\newcommand{\GF}[1]{\mbox{$\mbox{{\rm GF}}(#1)$}}
\newcommand{\abb}[3]{\mbox{$#1\,:\,#2\rightarrow#3$}}
\newcommand{\Abb}[5]{\mbox{$#1\,:\,#2\rightarrow#3,\;#4\mapsto #5$}}
\newcommand{\id}{\mbox{\rm id}}
\newcommand{\spn}{\mbox{{\rm span\,}}}
\newcommand{\Char}{\mbox{{\rm Char\,}}}

\newcommand{\vbf}{{\bf v}}

\newcommand{\Bcal}{{\cal B}}
\newcommand{\Ccal}{{\cal C}}
\newcommand{\Hcal}{{\cal H}}
\newcommand{\Ical}{{\cal I}}
\newcommand{\Kcal}{{\cal K}}
\newcommand{\Pcal}{{\cal P}}
\newcommand{\Rcal}{{\cal R}}
\newcommand{\Scal}{{\cal S}}
\newcommand{\Tcal}{{\cal T}}
\newcommand{\Vcal}{{\cal V}}
\newcommand{\Wcal}{{\cal W}}

\documentclass[a4paper]{article}
\usepackage{graphicx}
\usepackage{amssymb}
\begin{document}

%
%
%
%

\vspace*{3cm}
\noindent
{\Large{\bf Giuseppe Veronese and Ernst Witt -- Neighbours in \PG53}}

\abstand\noindent
{\sc Hans Havlicek}

\abstand\noindent
{\em Dedicated to J\'anos Acz\'el on the occasion of his 75th birthday}

   \abstand
   {
   \small\noindent
   {\bf Summary.}
   Let $P$ be a point of the Veronese surface $\Vcal$ in \PG53. Then there
   are four conics of $\Vcal$ through $P$. We show that the internal points
   of those conics form a $12$--cap which is a point model for Witt's
   $5$--$(12,6,1)$ design. In fact, this construction is ``dual'' to a
   similar construction that has been established in \cite{havl9x} recently.
   We give an explicit parametrization of the cap $\Kcal$; the domain is a
   dual affine plane which arises from \PG23 by removing one point. Thus, as
   a by--product, we obtain an easy approach to the extended ternary Golay
   code $G_{12}$. Finally, we discuss some other procedures that yield
   $12$--sets of points from the Veronese surface.
   }

\abstand\noindent
{\bf Mathematics Subject Classification (1991).} 51E22, 05B05.

\section{Introduction}

A construction of a cap $\Kcal$, in \PG53, which is a point model for {\em
Witt's $5$--$(12,6,1)$ design} $W_{12}$ (see, among others,
\cite[Chapter IV]{beth-jung-lenz85}) has been found by {\sc H.S.M.\ Coxeter}
\cite{coxe58} and, independently, by {\sc G.\ Pellegrino} \cite{pell74}. The
cap $\Kcal$ has exactly twelve points and any five distinct points of $\Kcal$
span a prime (hyperplane) of \PG53 which contains exactly six points
of $\Kcal$. Such a $\Kcal$ is projectively unique.
The group of collineations fixing $\Kcal$, as a set, is
the automorphism group of $W_{12}$, i.e.\ the {\em Mathieu group} $M_{12}$.
Also, {\sc J.A.\ Todd} \cite{todd59} has shown that there are exactly
twelve primes of \PG53 carrying no point of $\Kcal$. Those primes
gives rise to a point model $\Kcal^\ast$ of $W_{12}$ in the dual space of
\PG53.

The {\em Veronese surface} $\Vcal$ in \PG53 is a set of thirteen points;
cf.\ e.g.\ \cite[Chapter 25]{hirs-thas91}. It determines uniquely its {\em
dual Veronese surface} $\Vcal^\ast$ in the dual space of \PG53. As has been
pointed out in \cite{havl9x}, the following holds true: If one conic of the
Veronese surface $\Vcal$ is replaced with the set formed by the internal
points of that conic, then a point model $\Kcal$ of $W_{12}$ is obtained.
Figure \ref{bild} illustrates a conic in a projective plane of order three:
The four points of the conic, its three internal points, and its six external
points are drawn as squares, triangles, and hexagons, respectively.
%
{\unitlength1cm
   \begin{figure}[ht]\label{bild} 
      \begin{center}
         \begin {picture}(5.08,5.08)
         \put(0,0){\includegraphics[width=5.08cm]{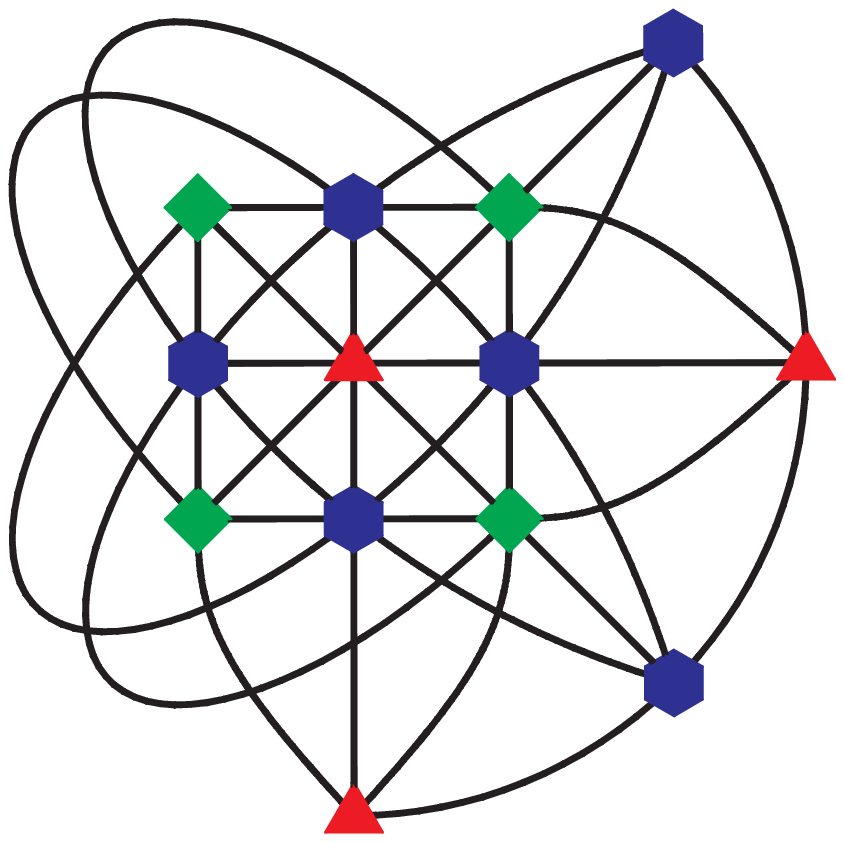}}
         \end{picture}
      \end{center}
      \caption{A conic in \PG23} 
   \end{figure}
}%

Clearly, now the question about the connection between the dual Veronese
surface $\Vcal^\ast$ and the dual point model $\Kcal^\ast$ arises. For our
purposes it will be convenient to adopt a dual point of view from the very
beginning: So we start with $\Vcal$, then go over to its dual $\Vcal^\ast$,
next apply the procedure of \cite{havl9x} to obtain a $\Kcal^\ast$ from
$\Vcal^\ast$, and finally go back to $\Kcal^{\ast\ast}=:\Kcal$.

In Theorem \ref{point-theo} we give a direct description (avoiding the dual
space) of that construction. It turns out that the present construction is
different from the one given in \cite{havl9x}, since now we do not have a
distinguished conic of $\Vcal$, but a distinguished point $P\in\Vcal$.
In contrast to \cite{havl9x}, $\Vcal$ and $\Kcal$ do not have common points,
but there is a bijection $\abb{\rho}{\Vcal\setminus\{P\}}\Kcal$ such that
$Y\in\Vcal\setminus\{P\}$, $Y^\rho\in\Kcal$, and $P$ are always collinear.

The {\em Veronese mapping} $\varphi$ is a bijection of \PG23 onto $\Vcal$. By
combining it with the above--mentioned bijection $\rho$, we find a bijection
of a dual affine plane (a \PG23 with one point removed) onto $\Kcal$. Using
homogeneous coordinates in \PG23 and \PG53 then gives an explicit parametric
representation $\psi$ of $\Kcal$ described in Theorem \ref{param-theo}. It is
a peculiar feature that one of the coordinate functions of $\psi$ contains an
{\em inhomogeneous\/} term, something which usually does not make sense. But
here it is meaningful, since $1$ is the only non--zero square in \GF3.

In Section \ref{replace} we discuss procedures, similar to that of Theorem
\ref{point-theo}, which yield $12$--sets of points from the Veronese surface
$\Vcal$. Some of them give point models of $W_{12}$, others do not. However,
one should keep in mind that each point model of $W_{12}$ arises in some way
or another from a fixed Veronese surface; so our discussion is far from being
complete. In fact, all our $12$--sets are ``close neighbours'' of the
Veronese surface $\Vcal$, as they belong to the algebraic hypersurface (of
order three) which is formed by all points of \PG53 that are on a chord of
$\Vcal$.

\section{The Veronese Surface and its Dual}\label{vero-dual}

Let us recall some properties  of  the  Veronese  surface  $\Vcal$  in  \PG53
\cite[Kapitel V]{bura61}, \cite{herz82}, \cite[Chapter 25]{hirs-thas91}:  The
term {\em conic plane\/} is used for a plane which meets $\Vcal$ in a conic.
Any two distinct conic planes have one and only one point in common. This
point belongs to $\Vcal$. Let $\Ccal$ be the set of all conics of $\Vcal$.
Then $(\Vcal,\Ccal,\in)$ is a \PG23.

For each $m\in\Ccal$ there is a unique prime which meets $\Vcal$ exactly in
$m$; this {\em osculating prime\/} (or {\em contact prime\/}) of $\Vcal$
along $m$ will be denoted by $\Hcal_m\Vcal$. A line $l$ is called a {\em
tangent\/} of $\Vcal$, if $l$ is a tangent of a conic $m\subset\Vcal$.
Given $P\in \Vcal$ the {\em tangent plane\/} of $\Vcal$ at $P$ is the union
of all tangents of $\Vcal$ which are running through $P$. It is written as
$\Tcal_P\Vcal$. Another description of osculating primes and tangent planes
is given by
   \begin{eqnarray}
   \label{vero-1}
   \Hcal_m\Vcal & = & \spn\left(\bigcup_{X\in m} \Tcal_X\Vcal\right)
   \quad (m\in\Ccal),\\
   \label{vero-2}
   \Tcal_P\Vcal & = & \bigcap_{P\in m\in\Ccal} \Hcal_m\Vcal
   \quad (P\in\Vcal).
   \end{eqnarray}
Any two distinct tangent planes have a unique common point; this point is not
on $\Vcal$.

If $\Scal$ is a subspace of \PG53, then let $[\Scal]^\ast$ be the star of
primes through $\Scal$. All osculating primes of $\Vcal$ form the {\em dual
Veronese surface} $\Vcal^\ast$, i.e.\ a Veronese surface in the dual space,
since $\Char\GF3\neq 2$. There is a one--one correspondence between the
tangent and conic planes of $\Vcal$ with the ``conic'' and ``tangent planes''
of $\Vcal^\ast$, respectively:

Each tangent plane $\Tcal_P\Vcal$ yields the ``conic plane''
$[\Tcal_P\Vcal]^\ast$ of $\Vcal^\ast$. The osculating primes of $\Vcal$ that
are passing through $P$ comprise the corresponding ``conic''
$c^\ast\subset\Vcal^\ast$. If we choose one ``point'' of the ``conic''
$c^\ast$, say $\Hcal_m\Vcal$, then its ``tangent'' is given by the pencil
$[\Tcal_P\Vcal\vee\spn m]^\ast$. An ``internal point'' of $c^\ast$ is a prime
$\Ical\in [\Tcal_P\Vcal]^\ast$ which is on no ``tangent'' of $c^\ast$, i.e.\
   \begin{equation}\label{vero-intern}
   \spn m \cap \Ical = \spn m \cap \Tcal_P\Vcal
   \mbox{ for all }m\in\Ccal \mbox{ with } P\in m.
   \end{equation}
Alternatively, an ``internal point'' of $c^\ast$ may be characterized as a
prime $\Ical$ of \PG53 satisfying
   \begin{equation} \label{vero-einpunkt}
   \Ical\cap\Vcal = \{P\},
   \end{equation}
since (\ref{vero-einpunkt}) implies that $\Ical$ corresponds to a quadric in
\PG32 consisting of one double point only, so that $\Tcal_P\Vcal\subset\Ical$
(cf.\ \cite[p.\ 168, Satz 1]{bura61}).

Likewise, each conic plane $\spn m$ ($m\in\Ccal$) yields the
``tangent plane'' $[\spn m]^\ast$ of $\Vcal^\ast$ at the ``point''
$\Hcal_m\Vcal\in\Vcal^\ast$.

\section{Point Models of $W_{12}$}

Let $\Kcal$ be a set of twelve points in a \PG53. Define a {\em block\/} of
$\Kcal$ as hyperplane section of $\Kcal$ which contains exactly six points.
Write $\Bcal$ for the set of all such blocks. If $(\Kcal,\Bcal,\in)$
is Witt's $5$--$(12,6,1)$ design $W_{12}$, then $\Kcal$ is called a {\em
point model\/} of $W_{12}$ in \PG53.

In what follows we put $\GF3=\{0,1,2\}=:F$.

   \begin{theo} \label{point-theo}
   Let $P$ be a point of the Veronese surface $\Vcal$ in \PG53. The four
   conics of $\Vcal$ through $P$ are denoted by $m_k$ ($k\in
   F\cup\{\infty\}$). The set of internal points of each $m_k$ is written as
   $\Delta_k$. Also let $c^\ast$ be the set of osculating primes of $\Vcal$
   through $P$ and $\Delta^\ast$ the set of all primes that meet $\Vcal$
   in $P$ only.
   Then the following holds true:
   \begin{enumerate}
   \item The set
      \begin{equation}\label{point-K}
      \Kcal := \bigcup_{k\in F\cup\{\infty\}} \Delta_k
      \end{equation}
   is a point model of the Witt design $W_{12}$.
   \item No point of $\Kcal$ is incident with a prime belonging to
      \begin{equation}
      \Kcal^\ast:=(\Vcal^\ast\setminus c^\ast)\cup\Delta^\ast.
      \end{equation}
   \end{enumerate}
   \end{theo}
\proof
According to Section \ref{vero-dual}, $c^\ast$ is a ``conic'' of the dual
Veronese surface and $\Delta^\ast$ is the set of its ``internal points''.
Hence $\Kcal^\ast$ is a point model of $W_{12}$ in the dual space of \PG53
\cite[Remark 3]{havl9x}. By a result of {\sc J.A.\ Todd} \cite[p.\
408]{todd59}, applied to $\Kcal^\ast$, there are exactly twelve points of
\PG53 which are not in any prime belonging to $\Kcal^\ast$. Moreover, those
points form a point model of $W_{12}$.

Each conic $m_k$ has exactly three internal points. The planes of two
distinct conics of the Veronese surface do not have common internal points.
Thus $\#\Kcal=12$ and the theorem follows, if we can show that no point of
$\Kcal$ lies in a prime belonging to $\Kcal^\ast$.

Let $\Ical$ be one of the three ``internal points'' of
$c^\ast$. By (\ref{vero-intern}), the prime $\Ical$ meets the plane of
each $m_k$ ($k\in F\cup\{\infty\}$) in the tangent of $m_k$ at $P$, so that
      \begin{equation}\label{point-schnittIK}
      \Ical\cap \Kcal = \emptyset.
      \end{equation}
Any of the remaining nine primes in $\Kcal^\ast$ is an osculating prime
$\Hcal_c\Vcal$ along a conic $c\in\Ccal\setminus\{m_0,m_1,m_2,m_\infty\}$. Put
$\{S_k\}:=c\cap m_k$. By (\ref{vero-2}),
$\Tcal_{S_k}\Vcal\subset\Hcal_c\Vcal$, whence $\Hcal_c\Vcal\cap \spn m_k $ is
the tangent of $m_k$ at $S_k$. It follows that
      \begin{equation}\label{point-schnittHK}
      \Hcal_c\Vcal\cap \Kcal = \emptyset,
      \end{equation}
which completes the proof.
\qed

\abstand\noindent
Now we introduce coordinates in order to obtain a parametric representation
of $\Kcal$. Assume that \PG23 and \PG53 are projective spaces $\Pcal(F^3)$
and $\Pcal(F^6)$, respectively. The Veronese mapping is given as
\begin{equation}
\label{vero}
\Abb\varphi{\Pcal(F^3)}{\Pcal(F^6)}
               {F(x_0,x_1,x_2)}{F(x_0^2,x_0 x_1,x_0 x_2,x_1^2,x_1
               x_2,x_2^2)}.
\end{equation}
We fix the points
   \begin{equation}\label{point-PU}
   U:=F(1,0,0) \mbox{ and } P:=U^\varphi = F(1,0,0,0,0,0).
   \end{equation}
Let $m_k\in\Ccal$ be a conic through $P$. Each bisecant of $m_k$ contains
exactly one internal point of $m_k$; see Figure \ref{bild}. Thus the mapping
   \begin{equation}\label{point-rho}
   \Abb\rho{\Vcal\setminus\{P\}}\Kcal{Y}{Y^\rho}
   \mbox{ with } P,Y,Y^\rho \mbox{ collinear }
   \end{equation}
is a well--defined bijection. Putting
   \begin{equation}\label{point-W}
   \Wcal:=\Pcal(F^3)\setminus\{U\}
   \end{equation}
yields the bijection
   \begin{equation}\label{point-psi}
   \Abb{\psi}{\Wcal}{\Kcal}{X}{X^{\varphi\rho}}
   \end{equation}
whose domain is a dual affine plane.

   \begin{theo} \label{param-theo}
   In terms of homogeneous coordinates the mapping (\ref{point-psi}) takes
   the form
      \begin{equation}\label{point-psikoo}
      F(x_0,x_1,x_2)\stackrel{\psi}{\mapsto}
      F(x_0^2+1,x_0 x_1,x_0 x_2,x_1^2,x_1 x_2,x_2^2)
      \end{equation}
   \end{theo}
\proof
At first we note that the use of the inhomogeneous term $x_0^2 + 1$ is not
ambiguous: In fact, if $\abb{q}{F^{n+1}}{F}$ is a quadratic form, then
   \begin{equation}
   (2x_0,\ldots,2x_n)^q =
   2^2\cdot(x_0,\ldots,x_n)^q =
   1\cdot(x_0,\ldots,x_n)^q.
   \end{equation}
Next choose a fixed pair $(x_1,x_2)\in F^2\setminus\{(0,0)\}$. The Veronese
image of the line joining $U$ with the point $F(0,x_1,x_2)\in\Pcal(F^3)$ is
a conic $m$ through $P$. Then $m$ comprises those four points which are
spanned by the vectors
   \begin{equation}
   \vbf_{u} := (u^2,u x_1,u x_2,x_1^2,x_1 x_2, x_2^2) \quad (u\in F),\quad
   \vbf_\infty := (1,0,0,0,0,0).
   \end{equation}
The three internal points of the conic $m$ are the three diagonal points of
the planar quadrangle $m$. We observe
   \begin{equation}
   \vbf_{u} + \vbf_\infty= 2\vbf_{u+1} + 2 \vbf_{u+2} \mbox{ for all }
   u\in F,
   \end{equation}
whence $F(u,x_1,x_2)^\psi = (F\vbf_{u})^\rho = F(\vbf_{u} + \vbf_\infty)$
for all $u\in F$.

As $(x_1,x_2)$ varies in $F^2\setminus\{(0,0)\}$, we obtain all four conics
$\{m_0,m_1,m_2,m_\infty\}$, since those conics can be relabelled in such a
way that $(x_1,x_2)$ yields the conic $m_k$ with $k=x_2/x_1$.
\qed

   \begin{remark}
   {\em
   The dual Veronese mapping $\varphi^\ast$ assigns to each line $l$ of \PG23
   with homogeneous coordinates $F(a_0,a_1,a_2)$ the prime of \PG53 with
   homogeneous coordinates
      \begin{equation}\label{point-dualvero}
      F(a_0^2, 2a_0a_1, 2a_0a_2, a_1^2, 2a_1a_2, a_2^2),
      \end{equation}
   since $l^{\varphi^\ast}$ equals the osculating prime of $\Vcal$ along the
   conic $l^\varphi$. The image of $\varphi^\ast$ is the dual Veronese
   surface $\Vcal^\ast$.
   The nine lines which are not running through $U$ are characterized by
   $a_0\neq 0$, whence their images under $\varphi^\ast$ are immediate from
   (\ref{point-dualvero}). By \cite[Remark 1]{havl9x}, the remaining three
   primes in $\Kcal^\ast$ have coordinates
   \begin{equation}\label{point-innenkoo}
   F(0,0,0,1,0,1),\;    F(0,0,0,2,2,1),\;    F(0,0,0,2,1,1).
   \end{equation}
   Note that the $01$--, $02$--, and $12$--coordinates in \cite{havl9x} have
   to multiplied by $2$ in order to match (\ref{point-dualvero}).
   By virtue of (\ref{point-psikoo}), (\ref{point-dualvero}), and
   (\ref{point-innenkoo}) it is easy to verify in terms of coordinates that
   no point of $\Kcal$ is incident with a prime belonging to $\Kcal^\ast$.
   This gives an alternative proof of Theorem \ref{point-theo}.
   }
   \end{remark}

   \begin{remark}
   {\em
   With the help of formula (\ref{point-psikoo}) one  may  immediately  write
   down twelve vectors of $F^6$ representing the points of $\Kcal$. If those
   vectors are regarded as columns of a $6\times 12$ matrix over $F$, then a
   generator matrix of the {\em extended ternary Golay code} $G_{12}$ arises.
   We refer to \cite{assm-matt66}, \cite[8.6]{assm-key92}, and \cite{havl9x}
   for further details on the connections between the Witt design $W_{12}$
   and coding theory.
   }
   \end{remark}

\section{Replacing Conics of the Veronese Surface} \label{replace}

In what follows we shall stick to the terminology introduced in Theorem
\ref{point-theo}, as we aim at generalizing the construction given there.

Choose one conic $m_k$ ($k\in F\cup\{\infty\}$). Let $t_k$
be the tangent of $m_k$ at $P$, $\Delta_{k,0}:=m_k\setminus \{P\}$, and
$\Delta_{k,1} := \Delta_k$. Denote by $\Delta_{k,2}$ the set of all
external points of $m_k$ that are off the tangent $t_k$. It is
easily seen from Figure \ref{bild} that there is a unique elation $\kappa_k$
of the plane $\spn m_k$ with centre $P$ and axis $t_k$ such that
   \begin{equation}
   (\Delta_{k,j})^{\kappa_k}=\Delta_{k,j+1} \mbox{ for all }j\in F.
   \end{equation}
Also, each $\Delta_{k,j+1}$ is the set of internal points of the conic
$\Delta_{k,j}\cup\{P\}$, whence the restrictions of $\kappa_k$ and $\rho$ to
the conic $m_k$ are coinciding.

All  four  collineations  $\kappa_k$  do  not  simultaneously  extend  to   a
collineation of \PG53, since  $\Vcal\setminus\{P\}$  contains  four  distinct
coplanar points, whereas  $\Kcal$  does  not.  However,  if  we  choose  three
collineations, e.g.\ $\kappa_1$, $\kappa_2$, and $\kappa_\infty$,  then  they
extend to a unique perspective collineation $\mu_0$ of \PG53 with centre  $P$
and axis $\Hcal_{m_0}\Vcal$: If the numbering of the  conics  $m_k$  is  done
according to the proof of Theorem \ref{param-theo}, then $0=0/1$ implies
that $\mu_0$ is given by
   \begin{equation}\label{123-erweiter}
   F(y_{00},y_{01},y_{02},y_{11},y_{12},y_{22})\mapsto
   F(y_{00}+y_{22},y_{01},y_{02},y_{11},y_{12},y_{22}).
   \end{equation}
Clearly, all points of the conic plane $\spn m_0\subset\Hcal_{m_0}\Vcal$ are
fixed under $\mu_0$, i.e.\ $\mu_0$ extends $\kappa_0^0=\id$.

One may assign to each $(p,q,r,s)\in F^4$ the quadruple
$(\kappa_0^p, \kappa_1^q, \kappa_2^r, \kappa_\infty^s$) of colline\-ations.
This mapping is an isomorphism of the group $(F^4,+)$ onto the direct sum of
the collineation groups generated by the $\kappa_k$'s. In addition, each
$(p,q,r,s)\in F^4$ yields the
point set
   \begin{equation}
   \Delta_{0,p}\cup\Delta_{1,q}\cup\Delta_{2,r}\cup\Delta_{\infty,s}
   \end{equation}
consisting of twelve points in \PG53.

Each permutation of the four conics $m_0,m_1,m_2,m_\infty$  arises  from   an
automorphic projective collineation of the Veronese surface $\Vcal$.  So,  if
two elements of $F^4$ differ only in the arrangement of their entries, then
they yield projectively equivalent $12$--sets.

The  subgroup  $S$  of  $(F^4,+)$  generated  by  $(1,1,1,0)$,   $(1,1,0,1)$,
$(1,0,1,1)$, and $(0,1,1,1)$ consists of all $(p,q,r,s)$ with $p+q+r+s=0$.
Each element of $S$ yields a quadruple of planar collineations that extend to
a collineation of \PG53 by (\ref{123-erweiter}). Hence all corresponding
$12$--sets are projectively equivalent to $\Vcal\setminus\{P\}$.

The elements of the coset $(1,0,0,0)+S$ are characterized by $p+q+r+s=1$ and
yield $12$--sets that are projectively equivalent to $\Kcal$.

We mention without proof some results about the $12$--sets that arise from
the elements of the remaining coset ($p+q+r+s=2$). Each such point set, say
$\Rcal$, is neither projectively equivalent to $\Kcal$ nor projectively
equivalent to a subset of $\Vcal$. Among the $364$ primes of \PG53 there are
exactly $42$ which meet $\Rcal$ in precisely six points. Those $42$ primes
have $P$ as their only common point. Thus $P$ is invariant under the group of
automorphic collineations of $\Rcal$. So it seems natural to project $\Rcal$
through the point $P$ to a prime $\Hcal$ of \PG53 not containing $P$: The
conic planes $m_k$ ($k\in F\cup\{\infty\}$) are projected to four mutually
skew lines $l_k\subset\Hcal$, the tangent plane $\Tcal_P\Vcal$ of the
Veronese surface goes over to the only transversal line of the $l_k$'s, and
the set $\Rcal$ is mapped onto those twelve points of the four lines $l_k$
that are not on their common transversal line.

\abstand\noindent
{\bf Acknowledgement}

\abstand\noindent
This research was supported by the Austrian National Science Fund (FWF),
project P12353--MAT. The author is obliged to Corrado Zanella for critical
remarks and warm hospitality at the `Universit\`a degli Studi di Padova'
(Italy) in March 1998.


\abstand
\noindent
Hans Havlicek\\
Abteilung f\"ur Lineare Algebra und Geometrie\\
Technische Universit\"at\\
Wiedner Hauptstra{\ss}e 8--10\\
A--1040 Wien\\
Austria\\
e-mail: {havlicek@geometrie.tuwien.ac.at}
\end{document}